\documentclass[leqno,12pt]{article}

\usepackage[usenames, dvipsnames]{xcolor}
\usepackage{tikz}

\usepackage{mathpazo}
\usepackage{amssymb,hyperref,amsthm}
\usepackage{euscript}
\usepackage{flafter}
\usepackage{pstricks}
\usepackage[latin1]{inputenc}
\usepackage{amsmath}
\usepackage{amsfonts}
\usepackage{pstricks-add}
\usepackage{yfonts}
\usepackage{egothic}
\usepackage{dsfont}
\usepackage{bm}
\usepackage{variations}

\hypersetup{
    linktoc=page,
   linkcolor=red,          
  citecolor=blue,        
    filecolor=blue,      
   urlcolor=cyan,
    colorlinks=true           
}

 \newtheorem{theorem}{Theorem}[section]

 \newtheorem{corollary}[theorem]{Corollary}

\usepackage[refpage]{nomencl}

\usepackage{makeidx}
\makeindex

\usepackage{mathrsfs} 

\evensidemargin\oddsidemargin

\newcommand{\eqnsection}{
\renewcommand{\theequation}{\thesection.\arabic{equation}}
   \makeatletter
   \csname  @addtoreset\endcsname{equation}{section}
   \makeatother}
\eqnsection

\def\P{{\bf P}}
\def\E{{\bf E}}

\def\T{{\mathbb T}}

\def\1{{\mathds{1}}}

\newcommand{\R}{\mathbb{R}}
\renewcommand{\T}{\mathbb{T}}

\newcommand{\ind}[1]{\mathbf{1}_{\left\{ #1 \right\}}}
\renewcommand{\epsilon}{\varepsilon}

\usepackage{tikz}

\title{\bf On the tail distribution of the solution to some law equation.}

\author{Xinxin CHEN\footnote{Institut Camille Jordan - C.N.R.S. UMR 5208 - Universit\'e Claude Bernard Lyon 1
(France).
 \newline \vspace{0.1cm} \hspace{0.2cm} $\dag$School of Mathematical Sciences and LPMC, Nankai University
(China). 
\newline \vspace{0.1cm}    MSC 2000  60E05 60G70. \newline \vspace{0.5cm} \textit{Key words :  distribution equation, multivariate regularly varying} } $\ $,$\ $ Chunhua Ma$^\dag$ }

\begin{document}

\maketitle
\begin{abstract}
We consider a distribution equation which was initially studied by Bertoin \cite{Bertoin}: 
\[
M \stackrel{d}{=} \max\{\widetilde{\nu}, \max_{1\leq k\leq \nu}M_k\}.
\] 
where $\{M_k\}_{k\geq 1}$ are i.i.d. copies of $M$ and independent of $(\widetilde{\nu}, \nu)\in\mathbb{R}_+\times\mathbb{N}$. We obtain the tail behaviour of the solution of a generalised equation in a different but direct method by considering the joint tail of $(\widetilde{\nu}, \nu)$.
\end{abstract}

\section{Introduction to questions}
\label{Intro}
For a random variable $\nu\in\mathbb{N}$, we consider the following equation of distributions on $\mathbb{N}$:
\begin{equation}\label{eq1}
M\stackrel{d}{=}\max\{\nu, \max_{1\leq k\leq \nu}M_k\}
\end{equation}
where $M_k$ are i.i.d. copies of $M\in\mathbb{N}$ and independent of $\nu$. In fact, for a Galton-Waston tree $\T$ with offspring $\nu$, let 
\[
M:=\sup_{u\in\T}\nu_u
\]
be the largest offspring. Clearly $M$ satisfies the equation \eqref{eq1}.

Note that if $\E[\nu]\leq 1$ then $M<\infty$ a.s. and that if $\E[\nu]>1$, $\P(M=\infty)=\P(T=\infty)$.

In this paper, we only consider the critical case when $\E[\nu]=1$. If $\E[\nu]<1$, the tail distribution of $M$ is of the same order as $\nu$. For two sequences $(a_n)_{n\geq1}$ and $(b_n)_{n\geq1}$, we write $a_n\sim b_n$ as $n\rightarrow \infty$ if $\lim_{n\rightarrow \infty}\frac{a_n}{b_n}=1$ and we write $a_n\asymp b_n$ if $K_1 b_n\leq a_n\leq K_2 b_n$ for some positive constants $K_1,K_2>0$.

Bertoin \cite{Bertoin} considered this equation \eqref{eq1} and by use of the link between critical Galton-Watson tree and centred random walk, he proved the following theorem.
\begin{theorem}[Bertoin]
\label{thm1}
For $M$ whose distribution satisfies the equation \eqref{eq1},
\begin{enumerate}
\item if $\P(\nu>n)=n^{-\alpha}\ell(n) $ with $\alpha\in(1,2)$ and $\ell(\cdot)$ a slowly varying function at infinity, then as $n\rightarrow\infty$,
\begin{equation}\label{eq1v1}
\P(M>n)\sim \frac{C_\alpha}{n},
\end{equation}
where $C_\alpha\in(0,\infty)$ is a constant which depends only on $\alpha$;
\item if $\sigma^2:=Var(\nu)<\infty$, then
\begin{equation}\label{eq1v2}
\P(M>n)\sim \sqrt{\frac{2}{\sigma^2}\P(\nu>n)}.
\end{equation}
\end{enumerate}
\end{theorem}

We will reprove this theorem by direct calculations, using the generating function of $\nu$. 

More generally, for a random vector $(\widetilde{\nu},\nu)$ which takes values in $\mathbb{R}_+\times\mathbb{N}$ such that $\E[\nu]=1$, let us consider the following equation of distribution:
\begin{equation}\label{eq2}
M \stackrel{d}{=} \max\{\widetilde{\nu}, \max_{1\leq k\leq \nu}M_k\}
\end{equation}
where $M_k$ are i.i.d. copies of $M\in\mathbb{R}_+$ and independent of $(\widetilde{\nu},\nu)$. The distribution of $M$ differs according to the joint distribution of $(\widetilde{\nu},\nu)$. We first consider some special cases in the following.

\begin{theorem}\label{thm1+}
For $M$ whose distribution satisfies the equation \eqref{eq2},
 if $\E[\nu]=1$ and $\E[\nu^2]<\infty$, then as $r\rightarrow\infty$,
\begin{equation}\label{eq2v2}
\P(M>r)\sim \sqrt{\frac{2}{\sigma^2}\P(\widetilde{\nu}>r)};
\end{equation}
where $\sigma^2=Var(\nu)$
\end{theorem}

For random vector $(\widetilde{\nu},\nu)$ where $\nu$ has infinite variance, we need to consider multivariate regularly varying condition. One can refer to Chapter 6 in \cite{Res09} for more details.

For a $d$-dimensional random vector $\mathbf{X}\in\mathbb{R}^d$, its law is regularly varying of index $\alpha\in(0,\infty)$ if for some norm $||\cdot||$ on $\mathbb{R}^d$, there exists a random vector $\theta$ on the unit sphere $\mathbb{S}^{d-1}=\{x\in\mathbb{R}^{d}\vert ||x||=1\}$ such that for any $u\in(0,\infty)$ and as $x\rightarrow\infty$,
\begin{equation}\label{regularvarying}
\frac{\P(||\mathbf{X}||>ux, \frac{\mathbf{X}}{||\mathbf{X}||}\in \cdot)}{\P(||\mathbf{X}||>x)}\xrightarrow{weak} u^{-\alpha}\P(\theta\in\cdot)
\end{equation}
where the convergence is on weak topology of finite measures, i.e. for $C_b(\mathbb{R}^d)$. 

The equivalent characterization of multivariate regular variations is as follows. Recall that a measurable function $V:(0,\infty)\rightarrow (0,\infty)$ is regularly varying of index $\rho\in\mathbb{R}$ if as $x\rightarrow\infty$,
\[
V(xy)/V(x)\rightarrow y^\rho, \forall y\in(0,\infty).
\]
A $d$-dimensional random vector $\mathbf{X}$ is then regularly varying of index $\alpha\in(0,\infty)$ if and only if there exists a regularly varying function $V$ of index $-\alpha$ and a nonzero Radon measure $\mu$ on $\mathbb{R}^d$ such that, as $x\rightarrow\infty$,
\begin{equation}\label{eq3}
\frac{1}{V(x)}\P(x^{-1}\mathbf{X}\in\cdot)\xrightarrow{vague} \mu(\cdot)
\end{equation}
where vague convergence is for all functions in $C_K^+(\mathbb{R}^d)$. Here the measure $\mu$ is homogeneous of order $-\alpha$. A choice for the function $V$ is that $V(x)=\P(||\mathbf{X}||>x)$ in which case the vague convergence in \eqref{eq3} is also weak convergence, and for all $u\in(0,\infty)$, by \eqref{regularvarying},
\[
\mu(\{x\in\mathbb{R}^d: ||x||>u, x/||x||\in \cdot\})=u^{-\alpha}\P(\theta\in\cdot).
\]
Then the restriction of $\mu$ on $\{x: ||x||>1\}$ is a probability measure. Hence, for any $f\in C_b(\mathbb{R}^d)$,
\[
\E[f(x^{-1}\mathbf{X})\Big\vert ||\mathbf{X}||>x]\xrightarrow{x\rightarrow\infty} \int f(y)1_{||y||>1}\mu(dy).
\]

Now we state the following theorem for jointly regularly varying setting.

\begin{theorem}\label{thm2}
For $M$ whose distribution satisfies the equation \eqref{eq2} and for $X=(\widetilde{\nu},\nu)$, suppose that $X$ is regularly varying of index $\alpha\in(1,2)$ associated with $V(x)=\P(||X||>x)$ and $\mu$ given in \eqref{eq3} such that
\[
\mu(\{y=(y_1,y_2)\in\mathbb{R}^2_+: y_1>0, y_2>1\})>0\textrm{ and }\mu(\{y=(y_1,y_2)\in\mathbb{R}^2_+: y_2=1<||y||\})=0.
\]
Suppose moreover that there exist $c_1, c_2\in(0,\infty)$ such that
\[
\frac{\P(\nu>x)}{\P(||X||>x)}\xrightarrow{x\rightarrow\infty} c_1, \textrm{ and }\frac{\P(\widetilde{\nu}>x)}{\P(||X||>x)}\xrightarrow{x\rightarrow\infty} c_2.
\]
Then, as $r\rightarrow\infty$,
\begin{equation}\label{eq2v1}
\P(M>r)\sim \frac{C_{\alpha,\mu}}{r}
\end{equation}
where $C_{\alpha,\mu}$ is the unique positive solution of $c_1\alpha\Gamma(-\alpha)x^\alpha-\int_{\mathbb{R}^2_+} e^{-xy_2}1_{y_1>1}\mu(dy)=0$. 
\end{theorem}

\textbf{Remark}: Instead of the vector $(\widetilde{\nu},\nu)$, if these assumptions hold for the vector $(\widetilde{\nu}^\gamma, \nu)$ for some $\gamma>0$, these arguments still work for $M^\gamma$. In this case, the distribution of $(\widetilde{\nu},\nu)$ is called non-standard regularly varying according to \cite{Res09}.

Let us state the following corollary by assuming the joint tail of $(\widetilde{\nu},\nu)$.

\begin{corollary}\label{cor}
For the random vector $(\widetilde{\nu},\nu)$ with $\E[\nu]=1$, assume that there exist a non-increasing function $\gamma:\mathbb{R}_+\rightarrow\mathbb{R}_+^*$ and $\alpha>0$ such that  for any $b\in\mathbb{R}_+$, as $r\rightarrow\infty$, 
\[
\P(\widetilde{\nu}\geq r, \nu\geq br)\sim \gamma(b) r^{-\alpha}, 
\]
and $\P(\nu\geq r)\sim c r^{-\alpha}$ with some $c\in(0,\infty)$. For random variable $M$ satisfying the equation \eqref{eq2},
\begin{enumerate}
\item if $\alpha\in(1,2)$, then
\begin{equation}\label{1st}
\P(M>r)\sim \frac{C_{\alpha,\gamma, c}}{r},\textrm{ as } r\rightarrow\infty.
\end{equation}
\item if $\alpha=2$, then
\begin{equation}\label{2nd}
\P(M>r)\sim \frac{C_{\alpha,\gamma, c}}{r\sqrt{\log r}},\textrm{ as } r\rightarrow\infty.
\end{equation}
\item if $\alpha>2$, then
\begin{equation}\label{3rd}
\P(M>r)\sim \sqrt{\frac{2}{\E[\nu^2]-1}\P(\widetilde{\nu}>r)},\textrm{ as } r\rightarrow\infty.
\end{equation}
\end{enumerate}

\end{corollary}

\textbf{Remark}: One can refer to \cite{CdR18} for an application and the motivation of this corollary.

Note that Theorem \ref{thm1} is a special case of Theorems \ref{thm2} and \ref{thm1+} by taking $\widetilde{\nu}=\nu$. We only need to prove the last two theorems. Let us explain the main idea of proof here, especially for Theorem \ref{thm2} when $(\widetilde{\nu},\nu)$ is regularly varying of index $\alpha\in(1,2)$. 

In fact, let $f(s):=\E[s^\nu]$ be the generating function of $\nu$. Observe from \eqref{eq2} that for any $r>0$,
\begin{align*}
\P(M\leq r)=&\P(\widetilde{\nu}\leq r; \max_{1\leq i\leq \nu}M_k\leq r)\\
=&\P\left(\max_{1\leq i\leq \nu}M_k\leq r\right)-\P\left(\widetilde{\nu}> r; \max_{1\leq i\leq \nu}M_k\leq r\right)\\
=&f\left(1-\P(M> r)\right)-\E\left[\left(1-\P(M>r)\right)^\nu; \widetilde{\nu}> r\right]
\end{align*}
where the last equality follows from the independence between $(\widetilde{\nu},\nu) $ and all $M_k$. Write $x_r:=\P(M>r)$ for convenience, we obtain that
\begin{equation}\label{eqxr}
1-x_r=f(1-x_r)-\E[(1-x_r)^\nu; \widetilde{\nu}>r].
\end{equation}
Inspired by this equation, we define for any $r>0$ and for any $x\in[0,r]$:
\[
\Phi_r(x):=[f(1-\frac{x}{r})-(1-\frac{x}{r})]-\E[(1-\frac{x}{r})^\nu; \widetilde{\nu}>r].
\]
Apparently, one sees that $\Phi_r(rx_r)=0$. On the one hand, we show the tightness of $rx_r$. On the other hand, for some positive deterministic sequence $(\gamma_r)$, we will show that $\gamma_r\Phi_r$ converges to some continuous function $\Phi$ uniformly on any compact of $[0,\infty)$ and that the limit function $\Phi$ has one unique zero in $[0,\infty)$. This implies that $rx_r$ converges to the unique zero of $\Phi$. 

If $\nu$ has finite variance, the proof will be similar by changing the rescaling term. 

The paper is organised as follows. In Section 2, we study the tail of $M$ when $\nu$ has finite variance and prove Theorem \ref{thm1+}. In Section 3, we prove Theorem \ref{thm2} and Corollary \ref{cor}.

We write $f(x)\asymp g(x)$ as $x\rightarrow x_0$ if $0<\liminf_{x\rightarrow x_0}\frac{f(x)}{g(x)}\leq \limsup_{x\rightarrow x_0}\frac{f(x)}{g(x)}<\infty$.

\section{Theorem \ref{thm1+}: \eqref{eq2} when $\E[\nu^2]<\infty$}\label{V2}

We begin with the simpler case when $\nu$ has finite second moment. Note that \eqref{eq1v2} is a special case of \eqref{eq2v2} by taking $\widetilde{\nu}=\nu$. So it suffices to prove \eqref{eq2v2} in Theorem \ref{thm1+}.

Recall that for $x_r=\P(M>r)$, we have the equation \eqref{eqxr}. We first show that $x_r\asymp\sqrt{\P(\widetilde{\nu}>r)}$.

For the upper bound, one sees that
\[
f(1-x_r)-(1-x_r)=\E[(1-x_r)^\nu; \widetilde{\nu}>r]\leq \P(\widetilde{\nu}>r).
\]
For the generating function $f$ of $\nu$ with $f'(1)=1$, let $b(s)$ be the function on $[0,1)$ such that
\begin{equation}\label{v2b}
b(s)=\frac{1}{1-f(s)}-\frac{1}{f'(1)(1-s)}=\frac{1}{1-f(s)}-\frac{1}{1-s}.
\end{equation}
Then as proved in Lemma 2.1 of \cite{GK00}, if $Var(\nu)<\infty$, one has
\begin{equation}\label{v2bb}
0\leq b(s)\leq \sigma^2, \forall 0\leq s<1, \textrm{ and } \lim_{s\uparrow1}b(s)=\sigma^2/2.
\end{equation}

Plugging $s=1-x_r$ in \eqref{v2b} then using \eqref{v2bb} gives us that
\[
f(1-x_r)-(1-x_r)=x_r-(\frac{1}{x_r}+b(1-x_r))^{-1}=\frac{x_r^2b(1-x_r)}{1+x_rb(1-x_r)}\geq \frac{b(1-x_r)}{1+\sigma^2} x_r^2.
\]
It follows that
\[
x_r^2\leq \P(\widetilde{\nu}>r) \frac{1+\sigma^2}{b(1-x_r)}.
\]
Note that $\P(M<\infty)=1$ as $\E[\nu]=1$. Apparently $x_r\downarrow 0$ as $r\uparrow\infty$. So, $\lim_{r\uparrow\infty}b(1-x_r)=\sigma^2/2$ by \eqref{v2b}. Then for $r$ sufficiently large, $b(1-x_r)\geq \sigma^2/4>0$. As a result, for $r$ large enough,
\[
x_r\leq 2\sqrt{\frac{1+\sigma^2}{\sigma^2}\P(\widetilde{\nu}>n)}.
\]

 For the lower bound, as $(1-x)^k\geq 1-kx$ for any $k\in\mathbb{N}$ and $x\in[0,1]$, we have
\begin{align}
f(1-x_r)-(1-x_r)=\E[(1-x_r)^\nu; \widetilde{\nu}>r]\geq& \E[(1-\nu x_r); \widetilde{\nu}>r]\nonumber\\
= &\P(\widetilde{\nu}>r)-x_r\E[\nu; \widetilde{\nu}>r].\label{lbd1}
\end{align}
Note that for the generating function $f$, we have $f'(1)=1$ and $f"(1)=\sigma^2$. So,
\[
f(1-x_r)-(1-x_r)=f(1-x_r)-f(1)+f'(1)x_r\leq \frac{f''(1)}{2}x_r^2.
\]
It follows that
\[
\frac{\sigma^2}{2}x_r^2\geq \P(\widetilde{\nu}>r)-x_r\E[\nu; \widetilde{\nu}>r]
\]
By Cauchy-Schwartz inequality, 
\begin{equation}\label{SCieq}
\E[\nu\ind{\widetilde{\nu}>r}]\leq \sqrt{\E[\nu^2\ind{\widetilde{\nu}>r}]\P(\widetilde{\nu}>r)}.
\end{equation}
Observe that $\E[\nu^2\ind{\widetilde{\nu}>r}]\rightarrow0$ as $r\rightarrow\infty$. Hence, $\E[\nu\ind{\widetilde{\nu}>r}]\leq\sqrt{\P(\widetilde{\nu}>r)}$ for $r$ large enough. This yields that
\[
\frac{\sigma^2}{2}x_r^2+x_r\sqrt{\P(\widetilde{\nu}>r)}\geq \P(\widetilde{\nu}>r).
\]
Consequently, for all sufficiently large $r$,
\[
x_r\geq \frac{-\sqrt{\P(\widetilde{\nu}>r)}+\sqrt{1+2\sigma^2}\sqrt{\P(\widetilde{\nu}>r)}}{\sigma^2}=\frac{2}{1+\sqrt{1+2\sigma^2}}\sqrt{\P(\widetilde{\nu}>r)}.
\]
Therefore, we obtain that
\[
x_r\asymp\sqrt{\P(\widetilde{\nu}>n)}.
\]
Next, let $\gamma_r:=\P(\widetilde{\nu}>r)$. We define for any $r>0$ and $x\in[0,1/\gamma_r]$,
\[
\Phi_r(x):=f\left(1-x\sqrt{\gamma_r}\right)-\left(1-x\sqrt{\gamma_r}\right)-\E\left[\left(1-x\sqrt{\gamma_r}\right)^\nu; \widetilde{\nu}>r\right]
\]
By \eqref{eqxr}, one has $\Phi_r(\frac{x_r}{\sqrt{\gamma_r}})=0$. We are going to show that 
\begin{equation}\label{unifcvg1}
\frac{1}{\gamma_r}\Phi_r(x)\xrightarrow{r\rightarrow\infty}\phi(x):= \frac{\sigma^2}{2}x^2-1,
\end{equation}
uniformly in any compact $K\subset[0,\infty)$. The pointwise convergence is trivial for $x=0$. 

We treat $f\left(1-x\sqrt{\gamma_r}\right)-\left(1-x\sqrt{\gamma_r}\right)$ and $\E\left[\left(1-x\sqrt{\gamma_r}\right)^\nu; \widetilde{\nu}>r\right]$ separately. First, for any $x\in(0,\infty)$, by \eqref{v2b} and \eqref{v2bb},
\[
n^2[f(1-\frac{x}{n})-(1-\frac{x}{n})]=\frac{b(1-\frac{x}{n})}{1+\frac{x}{n}b(1-\frac{x}{n})}x^2\rightarrow \frac{\sigma^2}{2}x^2,
\]
as $n\rightarrow\infty$. As $f(1-\frac{x}{n})-(1-\frac{x}{n})$ is monotone for $x$, Dini's theorem shows that this convergence is uniform in any compact $K\subset[0,\infty)$. Note that $\gamma_r\downarrow0$ as $r\uparrow\infty$. Replacing $n$ by  $1/\sqrt{\gamma_n}$, we get
\[
\frac{1}{\gamma_r}[f(1-x\sqrt{\gamma_r})-(1-x\sqrt{\gamma_r})]\rightarrow \frac{\sigma^2}{2}x^2.
\]
uniformly for $x$ in a compact $K\subset[0,\infty)$. It remains to show that for any $x\geq0$
\begin{equation}\label{rest}
\E\left[\left(1-x\sqrt{\gamma_r}\right)^\nu; \widetilde{\nu}>r\right]\sim \gamma_r.
\end{equation}
It is immediate that
\[
\E\left[\left(1-x\sqrt{\gamma_r}\right)^\nu; \widetilde{\nu}>r\right]\leq \P(\widetilde{\nu}>r)=\gamma_r.
\]
Similarly as \eqref{lbd1} and \eqref{SCieq}, 
\begin{align*}
\E\left[\left(1-x\sqrt{\gamma_r}\right)^\nu; \widetilde{\nu}>r\right]\geq& \P(\widetilde{\nu}>r)-x\sqrt{\gamma_r}\E[\nu; \widetilde{\nu}>r]\\
\geq &\gamma_r-x \gamma_r\sqrt{\E[\nu^2\ind{\widetilde{\nu}>r}]}.
\end{align*}
This implies that
\[
\Big\vert\E\left[\left(1-x\sqrt{\gamma_r}\right)^\nu; \widetilde{\nu}>r\right]-\gamma_r\Big\vert\leq x \gamma_r\sqrt{\E[\nu^2\ind{\widetilde{\nu}>r}]}.
\]
where $\lim_{r\rightarrow\infty}\gamma_r\sqrt{\E[\nu^2\ind{\widetilde{\nu}>r}]}=0$. Moreover, we also have the uniform convergence 
\[
\Big\vert\E\left[\left(1-x\sqrt{\gamma_r}\right)^\nu; \widetilde{\nu}>r\right]-\gamma_r\Big\vert\rightarrow 0
\]
 in any compact set on $\mathbb{R}_+$. We hence deduce that 
 \[
 \frac{1}{\gamma_r}\Phi_r(x)\rightarrow \phi(x).
 \]
uniformly in any compact set on $\mathbb{R}_+$.

Now let us prove the convergence of $\frac{x_r}{\sqrt{\gamma_r}}$ by showing  that any convergent subsequence converges towards the same limit. In fact, note that for any subsequence of $\{x_r\}$ such that as $k\rightarrow\infty$,
\[
\frac{x_{r_k}}{\sqrt{\gamma_{r_k}}}\rightarrow x^*\in\mathbb{R}_+, 
\]
the uniform convergence \eqref{unifcvg1} and the continuity of $\phi$ yield that
\[
\lim_{k\rightarrow\infty}\frac{1}{\gamma_{r_k}}\Phi_{r_k}\left(\frac{x_{r_k}}{\sqrt{\gamma_{r_k}}}\right)=\phi(x^*).
\]
By \eqref{eqxr},
\[
\phi(x^*)=0.
\]
So, $x^*=\sqrt{2/\sigma^2}$ and 
\[
x_r\sim \sqrt{\frac{2}{\sigma^2}\P(\widetilde{\nu}>r)}.
\]


\section{Theorem \ref{thm2}: \eqref{eq2} in the jointly regularly varying case}

In this section, we study the tail of $M$ given in the equation
\[
M\stackrel{d}{=} \max\{\widetilde{\nu}, \max_{1\leq k\leq \nu}M_k\},
\]
where $X=(\widetilde{\nu},\nu)$ is independent of  $M_k$, $k\geq1$, and has multivariate regularly varying tail. For $x_r=\P(M>r)$ with $r>0$, we first show that 
\[
x_r\asymp \frac{1}{r}, \textrm{ as } r\rightarrow\infty.
\]

According to the assumptions of Theorem \ref{thm2}, we can write $V(r)=r^{-\alpha}\ell(r)$ for some $\alpha\in(1,2)$ and $\ell$ slowly varying function, then 
\begin{equation}\label{hyp}
\P(\nu>r)=r^{-\alpha}\ell_1(r),\textrm{ and } \P(\widetilde{\nu}>r)=r^{-\alpha}\ell_2(r)
\end{equation}
where $\ell_1$ and $\ell_2$ are two slowly varying function at infinity such that $\ell_1(r)\sim c_1\ell(r)$ and $\ell_2(r)\sim c_2\ell(r)$ as $r\rightarrow\infty$.

We first show that $x_r=O( \frac{1}{r})$. In fact, by \eqref{eqxr}, one has
\begin{equation}\label{upxr}
x_r=1-f(1-x_r)+\E[(1-x_r)^\nu;\widetilde{\nu}>r]\leq 1-f(1-x_r)+\P(\widetilde{\nu}>r).
\end{equation}
Note that $f$ is the generating function of $\nu$, thus,
\begin{align*}
f(1-x_r)=\E[(1-x_r)^\nu]\geq &\E[(1-x_r)^\nu; \nu\leq r]\\
\geq &\E[(1-\nu x_r); \nu\leq r]\\
=& \P(\nu\leq r)-x_r\E[\nu; \nu\leq r].
\end{align*}
Plugging it into \eqref{upxr} yields that
\begin{align*}
x_r\leq &1-\P(\nu\leq r)+x_r\E[\nu; \nu\leq r]+\P(\widetilde{\nu}>r)\\
=&\P(\nu>r)+x_r\E[\nu; \nu\leq r]+\P(\widetilde{\nu}>r).
\end{align*}
Recall that $1=\E[\nu]$, so
\[
x_r(1-\E[\nu;\nu\leq r])=x_r\E[\nu; \nu>r]\leq \P(\nu>r)+\P(\widetilde{\nu}>r).
\]
By \eqref{hyp}, $\P(\widetilde{\nu}>r)\sim \frac{c_2}{c_1}\P(\nu>r)$. Moreover, the tail of $\nu$ is regularly varying of index $\alpha\in(1,2)$, then Karamata's theorem implies that $\E[\nu; \nu>r]\sim \frac{\alpha}{\alpha-1}\P(\nu>r)r$. Consequently,
$x_r=O(\frac{1}{r})$.

Apparently, $0\leq x_r r\leq C<\infty$ for some constant $C>0$. This means the tightness of $\{rx_r\}_{r>0}$. Now we turn to the uniform convergence of $\Phi_r$ defined as follows:
\begin{equation}\label{cvgpointwise}
\Phi_r(x)=f(1-\frac{x}{r})-(1-\frac{x}{r})-\E[(1-\frac{x}{r})^\nu; \widetilde{\nu}>r], \forall r> x\geq0.
\end{equation}
Clearly, $\Phi_r(0)\sim -c_2 r^{-\alpha}\ell(r)$. 

On the one hand, by Bingham {\it et al\,.}\cite[Theorem 8.1.6; P 333]{BGT87}, as $\P(\nu>r)=r^{-\alpha}\ell_1(r)$, we get that for $x>0$, 
\begin{equation}\label{cvgpointwise1}
f(1-\frac{x}{r})-(1-\frac{x}{r})\sim\alpha\Gamma(-\alpha)\Big(\frac{x}{r}\Big)^\alpha \ell_1\Big(\frac{r}{x}\Big)\sim c_1\alpha\Gamma(-\alpha) x^\alpha V(r), \quad r\rightarrow\infty,
\end{equation}
On the other hand, we study the convergence of $\frac{1}{V(r)}\E[(1-\frac{x}{r})^\nu; \widetilde{\nu}>r]$ by use of the multivariate regularly varying tail of $X=(\widetilde{\nu},\nu)$. Note that 
\begin{align*}
\frac{1}{V(r)}\E[(1-\frac{x}{r})^\nu; \widetilde{\nu}>r]=&\E[(1-\frac{x}{r})^\nu; \widetilde{\nu}>r\vert ||X||>r]\\
\leq & \E[e^{-x \frac{\nu}{r}}1_{\{\widetilde{\nu}>r\}}\vert ||X||>r]
\end{align*}
For the lower bound, note that for any $\epsilon\in(0,1)$ and $x>0$, there exists $r(x,\epsilon)>0$ such that for all $r\geq r(x,\epsilon)$, one has $(1-\frac{x}{r})\geq e^{-(1+\epsilon)\frac{x}{r}}$. Therefore,
\[
\frac{1}{V(r)}\E[(1-\frac{x}{r})^\nu; \widetilde{\nu}>r]\geq \E[e^{-(1+\epsilon)x \frac{\nu}{r}}1_{\{\widetilde{\nu}>r\}}\vert ||X||>r].
\]
The multivariate regularly tail of the random vector $X$ implies the weak convergence of $\P((\frac{\widetilde{\nu}}{r}, \frac{\nu}{r})\in\cdot\vert ||X||>r)$ towards $\mu_+(\cdot):=\mu(\cdot\cap\{y=(y_1,y_2)\in\mathbb{R}^2_+: ||y||>1\})$. Moreover, for $A:=(1,\infty)\times\mathbb{R}_+$, $\mu_+(\partial A)=0$ by the assumption. As a consequence, for any $x>0$,
\begin{align*}
\lim_{r\rightarrow\infty}\frac{\P(\widetilde{\nu}>r)}{V(r)}=& \mu_+(A)=c_2\\
\liminf_{r\rightarrow\infty}\frac{1}{V(r)}\E[(1-\frac{x}{r})^\nu; \widetilde{\nu}>r]\geq &\int e^{-(1+\epsilon )x y_2}1_{\{y_1>1\}}\mu_+(dy)\\
\limsup_{r\rightarrow\infty}\frac{1}{V(r)}\E[(1-\frac{x}{r})^\nu; \widetilde{\nu}>r]\leq &\int e^{-x y_2}1_{\{y_1>1\}}\mu_+(dy)
\end{align*}
Letting $\epsilon>0$ yields that
\begin{equation}\label{cvgpointwise2}
\lim_{r\rightarrow\infty}\frac{1}{V(r)}\E[(1-\frac{x}{r})^\nu; \widetilde{\nu}>r]=\int e^{-x y_2}1_{\{y_1>1\}}\mu_+(dy).
\end{equation}
Note that $f(1-\frac{x}{r})-(1-\frac{x}{r})$ is non-decreasing in $x$ and $\frac{1}{V(r)}\E[(1-\frac{x}{r})^\nu; \widetilde{\nu}>r]$ is non-increasing in $x$. So, by Dini's theorem, both the convergences \eqref{cvgpointwise1} and \eqref{cvgpointwise2} are uniform on any compact in $[0,\infty)$. Going back to \eqref{cvgpointwise}, one obtains that 
\[
\frac{1}{V(r)}\Phi_r(x)\xrightarrow{r\rightarrow\infty}\Phi(x):=c_1\alpha\Gamma(-\alpha) x^\alpha - \int e^{-x y_1}1_{\{y_2>1\}}\mu_+(dy),
\]
uniformly on any compact in $[0,\infty)$. Note that by dominated convergence theorem, $\Phi$ is continuous for $x\geq0$. Moreover, $\Phi$ is strict increasing on $\R_+$. The zero $C_{\alpha,\mu}>0$ of $\Phi$ exists and is unique because 
\[
\Phi(0)=-c_2<0\textrm{ and }\lim_{x\rightarrow\infty}\Phi(x)=\infty.
\]
The equation \eqref{eqxr} means that $\Phi_r(rx_r)=0$. So we conclude that $rx_r\rightarrow C_{\alpha,\mu}$.

\section{Proof of Corollary \ref{cor}}

It is clear that when $\alpha>2$, \eqref{3rd} follows from Theorem \ref{thm1+} and that when $\alpha\in(1,2)$, \eqref{1st} follows from Theorem \ref{thm2}. We only need to prove \eqref{2nd} for $\alpha=2$. 

First note that $\P(v\geq r)\sim c r^{-\alpha}$ and $\P(\widetilde{\nu}\geq r)\sim \gamma(0) r^{-\alpha}$. By Theorem 1.5 in \cite{CdRH16}, one has for $r$ sufficiently large,
\[
p_r:=\P(M\geq r)\asymp \frac{1}{r\sqrt{\log r}}.
\]
For $r\geq e$ and $x\in[0, r\sqrt{\log r}]$, define
\[
\Phi_r(x)=f(1-\frac{x}{r\sqrt{\log r}})-(1-\frac{x}{r\sqrt{\log r}})-\E\left[\left(1-\frac{x}{r\sqrt{\log r}}\right)^\nu; \widetilde{\nu}\geq r\right].
\]
Similarly, one has for any $x\geq0$, as $r\rightarrow\infty$,
\[
f(1-\frac{x}{r\sqrt{\log r}})-(1-\frac{x}{r\sqrt{\log r}})\sim C_2\left(\frac{x}{r\sqrt{\log r}}\right)^2\log\left(\frac{r\sqrt{\log r}}{x}\right)\sim C_2\frac{x^2}{r^2}.
\]
For any $\varepsilon\in(0,1)$ and $r\gg 1$, 
\[
\E\left[e^{-\varepsilon \frac{x\nu}{r}};\widetilde{\nu}\geq r\right]\leq \E\left[\left(1-\frac{x}{r\sqrt{\log r}}\right)^\nu; \widetilde{\nu}\geq r\right]\leq \P(\widetilde{\nu}\geq r)
\]
Apparently $r^2\P(\widetilde{\nu}\geq r)\rightarrow c$ as $r\rightarrow\infty$. For the lower bound, we see that for any $x>0$,
\begin{align*}
r^2\E\left[1-e^{-\varepsilon \frac{x\nu}{r}};\widetilde{\nu}\geq r\right]=&\int_0^\infty e^{-u} r^2\P(\widetilde{\nu}\geq r, \nu\geq r \frac{u}{\varepsilon x})du\\
\leq & \int_0^\infty e^{-u} du \left(r^2\P(\widetilde{\nu}\geq r)\right)
\end{align*}
As $\P(\widetilde{\nu}\geq r)\sim c r^{-\alpha}$ and $\P(\widetilde{\nu}\geq r,\nu\geq br)\sim \gamma(b) r^{-\alpha}$, by dominated convergence theorem, 
\[
\lim_{r\rightarrow\infty}r^2\E\left[1-e^{-\varepsilon \frac{x\nu}{r}};\widetilde{\nu}\geq r\right]=\int_0^\infty \gamma(\frac{u}{\varepsilon x})e^{-u}du
\]
Here one can show that $\gamma(b)\downarrow 0$ as $b\uparrow \infty$. In fact, for any $b>0$, we have
\[
r^2\P(\widetilde{\nu}\geq r,\nu\geq br)\leq r^2\P(\nu \geq br).
\]
Taking limit shows that $\gamma(b)\leq c b^{-\alpha}$. Consequently,
\[
\lim_{r\rightarrow\infty}r^2\E\left[e^{-\varepsilon \frac{x\nu}{r}};\widetilde{\nu}\geq r\right]=\int_0^\infty[\gamma(0)-\gamma(\frac{u}{\varepsilon x})]e^{-u}du\geq \int_0^\infty [\gamma(0)- \frac{c\varepsilon^2x^2}{u^2}\wedge \gamma(0)]e^{-u}du.
\]
Note that by dominated convergence theorem,
\[
\int_0^\infty \frac{c\varepsilon^2x^2}{u^2}\wedge \gamma(0)e^{-u}du\rightarrow 0\textrm{ as } \varepsilon\rightarrow 0.
\]
Therefore,
\begin{multline*}
\gamma(0)+o_\varepsilon(1)\leq \liminf_{r\rightarrow\infty}r^2\E\left[\left(1-\frac{x}{r\sqrt{\log r}}\right)^\nu; \widetilde{\nu}\geq r\right]\\
\leq \limsup_{r\rightarrow\infty}r^2\E\left[\left(1-\frac{x}{r\sqrt{\log r}}\right)^\nu; \widetilde{\nu}\geq r\right]\leq \gamma(0)
\end{multline*}
We end up with
\[
r^2\Phi_r(x)\rightarrow \phi(x):=C_2 x^2-\gamma(0).
\]
One can check the uniform convergence in any compact set on $\mathbb{R}_+$ by Dini's theorem. We hence conclude that $r\sqrt{\log r}p_r$ converges to $\sqrt{\frac{\gamma(0)}{C_2}}$.



\end{document}